\documentclass[12pt]{amsart}    
\usepackage{amscd}

\def\cal#1{\mathcal{#1}}

\def\NZQ{\Bbb}               

\def\ZZ{{\NZQ Z}}

\def\frk{\frak}               

\def\mm{{\frk m}}

\def\opn#1#2{\def#1{\operatorname{#2}}} 
\opn\chara{char}
\opn\length{\ell}
\opn\pd{pd}
\opn\rk{rk}
\opn\projdim{proj\,dim}
\opn\rank{rank}
\opn\depth{depth}
\opn\grade{grade}
\opn\height{height}
\opn\embdim{emb\,dim}
\opn\codim{codim}

\opn\Tr{Tr}
\opn\bigrank{big\,rank}
\opn\superheight{superheight}\opn\lcm{lcm}
\opn\trdeg{tr\,deg}%
\opn\reg{reg}
\opn\lreg{lreg}
\opn\div{div}
\opn\Div{Div}
\opn\cl{cl}
\opn\Cl{Cl}
\opn\Spec{Spec}
\opn\Supp{Supp}
\opn\supp{supp}
\opn\Sing{Sing}
\opn\Ass{Ass}
\opn\Ann{Ann}
\opn\Rad{Rad}
\opn\Soc{Soc}
\opn\Ker{Ker}
\opn\Coker{Coker}
\opn\Im{Im}
\opn\Hom{Hom}
\opn\Tor{Tor}
\opn\Ext{Ext}
\opn\End{End}
\opn\Aut{Aut}
\opn\id{id}

\opn\nat{nat}
\opn\pff{pf}
\opn\Pf{Pf}
\opn\GL{GL}
\opn\SL{SL}
\opn\mod{mod}
\opn\ord{ord}
\opn\aff{aff}
\opn\con{conv}
\opn\relint{relint}
\opn\st{st}
\opn\lk{lk}
\opn\cn{cn}
\opn\core{core}
\opn\vol{vol}
\opn\link{link}
\opn\star{star}
\opn\gr{gr}

\def\poly#1#2#3{#1[#2_1,\dots,#2_{#3}]}
\def\pot#1#2{#1[\kern-0.28ex[#2]\kern-0.28ex]}

\opn\dirlim{\underrightarrow{\lim}}
\opn\inivlim{\underleftarrow{\lim}}
\let\dirsum=\oplus

\let\iso=\cong

\let\Dirsum=\bigoplus

\let\mcone= * 
\let\to=\rightarrow
\def\namedTo#1{\overset{#1}{\longrightarrow}}
\let\To=\longrightarrow
\def\Implies{\ifmmode\Longrightarrow \else
     \unskip${}\Longrightarrow{}$\ignorespaces\fi}
\def\implies{\ifmmode\Rightarrow \else
     \unskip${}\Rightarrow{}$\ignorespaces\fi}
\def\iff{\ifmmode\Longleftrightarrow \else
     \unskip${}\Longleftrightarrow{}$\ignorespaces\fi}

\let\:=\colon
\newtheorem{Theorem}{Theorem}
\newtheorem{Lemma}[Theorem]{Lemma}
\newtheorem{Corollary}[Theorem]{Corollary}
\newtheorem{Proposition}[Theorem]{Proposition}
\newtheorem{Remark}[Theorem]{Remark}

\newtheorem{Example}[Theorem]{Example}

\newtheorem{Definition}[Theorem]{Definition}

\let\epsilon\varepsilon
\let\kappa=\varkappa
\opn\inii{in}
\opn\inim{inm}
\opn\set{set}
\def\pnt{{\raise0.5mm\hbox{\large\bf.}}}

\begin{document}

\title{b-sequences and approximations of 
generalized Cohen-Macaulay ideals}
\author{Yukihide Takayama}
\address{Yukihide Takayama, Department of Mathematical
Sciences, Ritsumeikan University, 
1-1-1 Nojihigashi, Kusatsu, Shiga 525-8577, Japan}
\email{takayama@se.ritsumei.ac.jp}
\date{}

\def\Coh#1#2{H_{\mm}^{#1}(#2)}

\newcommand{\AppTh}{Theorem~\ref{approxtheorem} }
\def\da{\downarrow}
\newcommand{\ua}{\uparrow}
\newcommand{\namedto}[1]{\buildrel\mbox{$#1$}\over\rightarrow}
\newcommand{\bdel}{\bar\partial}
\newcommand{\proj}{{\rm proj.}}

\maketitle

\begin{abstract} We introduce the notion of b-sequence for finitely
generated modules over Noetherian rings, which characterizes long
Bourbaki sequences. Our main concern is an application of this notion
to generalized Cohen-Macaulay approximation, which we introduced in
\cite{HT2}. We will show how we can construct long Bourbaki sequences
of non-trivial type characterizing generalized Cohen-Macaulay rings
by finding  suitable b-sequences.
\end{abstract}

\section*{Introduction}

Relation between Bourbaki sequences and local cohomomogies
has been studied several times, for example, by Evans-Griffith \cite{EG}
and Auslander-Buchweitz \cite{AB}. In \cite{HT2}, we studied  approximations
of generalized Cohen-Macaulay modules by (non-CM) maximal generalized 
Cohen-Macaulay modules.

Let $(R,\mm)$ be a Gorenstein local ring and consider a generalized
Cohen-Macaulay ideal $I\subset R$ of codimension~$r (\geq 2)$, which
is an ideal such that the local cohomology is 
$\Coh{i}{R/I} \iso M_i$ $i=0,\ldots, n-r-1$, for
some finite length $R$-modules $M_i$. Then there exists 
a maximal generalized Cohen-Macaulay module $M$ fitting into a 
length $r$ long Bourbaki sequence
\begin{equation}
\label{appseq}
  0 \To F_{r-1}\To\cdots\To F_1 \To M \To I \To 0 \quad (exact)
\end{equation}
where $F_i$ are $R$-free modules, such that 
$\Coh{i}{M} \iso \Coh{i}{I}$ for $i\leq n-r$ and $\Coh{n-r+1}{M} = 0$.
If we restrict ourselves to consider
$M$ satisfying the additional homological condition
\begin{equation}
\label{uniquenesscond}
   \Coh{i}{M} = 0 \qquad (n-r+2\leq i\leq n-1),
\end{equation}
then $M$ is unique up to direct $R$-free summands 
(See Proposition~\ref{ama-her-tak}).
We will call the Bourbaki sequence $(\ref{appseq})$ an
{\em approximation sequence} and $M$ an {\em approximation
module} of $I$.
The proof of this fact is carried out by constructing, with a
homological method, approximation modules $M$ that always satisfy the
condition (\ref{uniquenesscond}). 

In this paper, we are interested in approximation sequences
that do not satisfy the condition (\ref{uniquenesscond}).
We do not know a systematic method to construct such sequences,
particularly in the case of $r\geq 3$.
Recall that length $2$ Bourbaki sequence can be constructed
by finding basic elements (\cite{Bour} Chapter~VII \S 4). 
In section~1 we introduce the notion of {\em b-sequences} 
for Bourbaki sequences of arbitrary length, which plays 
a similar role to basic elements. Then we give a characterization 
of long Bourbaki sequences in terms of b-sequences (Theorem~\ref{maincor}).
Section~2 gives a characterization of (non-trivial) approximation
sequences that do not satisfy the condition (\ref{uniquenesscond}) in a
typical case in terms of 
b-sequences (Theorem~\ref{theorem:non-trivialcond}).
Some examples in the case of $r=3$ are considered in
section~3, where we focus on the special case of approximation
modules $M$ such that $\Coh{t+1}{M} = \Coh{n-1}{M} = K$ (field)
and $\Coh{i}{M} = 0$ otherwise for $i<n$. 

For a set $S$, we will denote by $\langle S\rangle$ the module generated by $S$.
Also, for a module $M$ over a ring $R$,  the  $i$th syzygy module will be 
denoted by $\Omega_i(M)$.   

\section{b-Sequences for Modules}
\def\bseq{b-sequence\;}

\subsection{b-sequences and long Bourbaki sequences}
Recall that length~2 Bourbaki sequences over a normal domain $R$
\begin{equation*}
	0 \To F \To M \To I \To 0,
\end{equation*}
where $F$ is a $R$-free module, $M$ is a finitely generated
torsion-free $R$ module and $I \subset R$ is an ideal, 
can be constructed by finding basic elements in $M$ 
(\cite{Bour} Chapter~VII \S 4). 
In this section, we introduce the notion of {\em \bseq},
which is a couterpart of basic elements  for long Bourbaki sequences
\begin{equation*}
   0\To F_{r-1}\To \cdots \To F_1 \To M \To I,
\end{equation*}
in particular for $r\geq 3$.
We first prove 
\def\keylemma{Lemma~\ref{keylemma}\;}
\begin{Lemma}
\label{keylemma}
Let $R$ be a Noetherian  ring and $M$ be 
a finitely generated $R$-module with a  presentation
$0\To\Ker\epsilon\To U\namedTo{\epsilon} M\To 0$,
with a finitely generated free $R$-module $U$.
Also let $f: F\To G$ be a monomorphism of $R$-modules
where $G$ is free of $\rank G = q$.
Then, following are equivalent.
\begin{enumerate}
\item [$(i)$] We have an exact  sequence
\[
 0 \To F
   \namedTo{f} G
   \To  M
   \To  I
   \To 0 \quad (exact)
\]
for an ideal $I \subset R$.
\item [$(ii)$]
We have 
$\beta_1,\ldots,\beta_q \in U \backslash \Ker\epsilon$
and $\varphi \in \Hom_R(U, R)$ such that 
\begin{enumerate}
\item [$(a)$] 
      $\Ker(\varphi)
	= \langle\beta_1,\ldots\beta_q\rangle + \Ker\epsilon$, and 
\item [$(b)$] we have the following commutative diagram
\[
\begin{CD}
     @.    0              @.            0         \\
  @.        @VVV                         @VVV      \\
  0  @>>>  \Ker \beta\circ f             @>>>  \Ker\beta  @>>> 0\\
  @.        @VVV                         @VVV             \\
  0  @>>>  F                             @>{f}>>  G                \\
  @.       @V{\beta\circ f}VV            @V{\beta}VV \\
0 @>>> \langle\beta_1,\ldots,\beta_q\rangle
      \cap \Ker\epsilon
                               @>>>  \langle\beta_1,\ldots,\beta_q\rangle \\
 @.        @VVV                         @VVV \\
     @.    0            @.               0   \\
\end{CD}
\]
where $\beta(m_i) = \beta_i$ for all $i$ with $\{m_1,\ldots, m_q\}$ a free basis 
of $G$.
\end{enumerate}
In this case we have $I = \Im\varphi$.
\end{enumerate}
\end{Lemma}
\begin{proof}
We first prove $(ii)$ to $(i)$. 
We set $I = \varphi(U)$. Then by $(a)$ we have the diagram
\begin{equation*}
\label{construction:proof-1}
\begin{CD}
0 @>>> \langle\beta_1,\ldots,\beta_q\rangle
       + \Ker\epsilon              @>>> U                         @>{\varphi}>> I @>>> 0 \\
@.     @V{\epsilon}VV              @V{\epsilon}VV                            \\
0 @>>> \langle\epsilon(\beta_1),\ldots,\epsilon(\beta_q)\rangle
                              @>>> M                           \\
@.      @VVV                        @VVV                                @.  \\
   @.    0                    @.     0                       @.            
\end{CD}
\end{equation*}
Then we can define a well-defined 
map $\psi : M \to I$ by $\psi(\epsilon(a)) = \varphi(a)$ for
all $a\in U$, and we have 
$\Ker\psi = \langle\epsilon(\beta_1),\ldots,\epsilon(\beta_q)\rangle$.
Also we define $g : G \To M$ by $g = \epsilon\circ\beta$.
Then, from the above diagram, we have an exact sequence
\begin{equation*}
	G \namedTo{g} M \namedTo{\psi} I \To 0.
\end{equation*}
On the other hand, we have 
\begin{eqnarray*}
 \Ker g   & = & 
   \{ \sum_{i=1}^q h_i m_i \mid \epsilon(\sum_{i=1}^q h_i\beta_i) = 0, h_i\in R \} 
\\                   & = &
   \{ \sum_{i=1}^q h_i m_i \mid \sum_{i=1}^q h_i\beta_i\in\Ker\epsilon, h_i \in R\}
\\                   & = & 
   \{ \sum_{i=1}^q h_i m_i \mid \beta(\sum_{i=1}^q h_im_i)
                                 \in\langle\beta_1,\ldots, \beta_q\rangle
				\cap  \Ker\epsilon\}
\\                   & = & 
   \{ u\in G  \mid \beta(u) \in\langle\beta_1,\ldots, \beta_q\rangle\cap \Ker\epsilon\}
\end{eqnarray*}
Now let $u\in G$ be such that 
$\beta(u)\in\langle\beta_1,\ldots,\beta_q\rangle\cap\Ker\epsilon$. 
Then $u$ must be in $f(F)$. In fact, by (b) we can choose 
$v\in F$ such that
$(\beta\circ f)(v) = \beta(u)$. Thus $u- f(v) \in\Ker\beta 
\iso \Ker (\beta\circ f)\subset F$, and we have $u \in f(v) + \Ker \beta 
\subset f(F)$ as required.
Thus we have $\Ker g \subset \Im f$ and the converse inclusion is clear
by (b).
Consequently, we have a desired exact  sequence.

Next we prove $(i)$ to $(ii)$. Given an exact  sequence
\[
  0 \To F \namedTo{f} G \namedTo{g} M \namedTo{\psi} I \To 0
\]
with an ideal $I\subset R$. 
Then $\psi \in \Hom_R(M, R)$ and we have
$\Ker\epsilon + N = \epsilon^{-1}(\Ker\psi) (\subset U)$ for some submodule $N(\ne 0)$ of
$U$. Then we can choose a finite set of generators $\{\beta_i\}_i$ of $N$ such that 
$\epsilon(\beta_i) = g(m_i)$ $(\forall i)$ where $\{m_i\}_i$ is a 
$R$-free basis of $G$.
Then we have the following diagram:
\[
\begin{CD}
   @.                                           @.      0        @. \\
@.       @.                                             @AAA     @. \\
   @.    G                                      @>{g}>> M        @>{\psi}>> I    @>>> 0 \\
@.       @V{\beta}VV                                  @A{\epsilon}AA  \\ 
0 @>>> \langle\beta_1,\ldots,\beta_q\rangle @>>> U   \\  
@.       @VVV                                         @.               @. \\
    @.     0                                    @.           @.           \\
\end{CD}
\]
where we define $\beta(m_i) = \beta_i$ $(\forall i)$. 
Then by defining $\varphi = \psi\circ \epsilon$, we have 
$\{\beta_i\}_i$ and $\varphi \in \Hom_R(U,R)$ satisfying the condition $(ii)(a)$.
Now we prove (ii)(b).
Since  $\Ker g =\Ker (\epsilon\circ\beta) = \Im f \iso F$ 
we readily have the following diagram:
\[
\begin{CD}
   @.    0               @.       0          @.    \\
@.      @VVV                      @VVV           @.\\
0  @>>> \Ker\beta\circ f @>>>     \Ker\beta @>>> 0 \\
@.      @VVV                      @VVV           @.\\
0 @>>> F                 @>{f}>> G @>{g}>> M\\
@.       @.                       @V{\beta}VV     @| \\
0 @>>> \langle\beta_1,\ldots,\beta_q\rangle\cap\Ker\epsilon
                         @>>>     \langle\beta_1,\ldots,\beta_q\rangle
                                                  @>{\epsilon\mid_{\Im\beta}}>>
                                                  M\\
@.       @.                       @VVV             @. \\
    @.                   @.       0               @. \\
\end{CD}
\]
Notice that since $\Ker\beta \subset\Ker (\epsilon\circ\beta) 
= \Ker g = \Im f$ we have the exactness of the first row.
Since 
$\Im (\beta\circ f) = \beta (\Ker(\epsilon\circ \beta)) 
=  \Ker (\epsilon\mid_{\Im\beta})
= \langle\beta_1,\ldots,\beta_q\rangle\cap\Ker\epsilon$, 
we have a well-defined surjection
\[
   \beta\circ f :  F  \To \langle\beta_1,\ldots,\beta_q\rangle\cap\Ker\epsilon
\]
as required.
\end{proof}

Now we introduce the notion of \bseq.
\begin{Definition}
Let $R$ be a Noetherian ring.
For a finitely generated $R$-module $M$ with a  presentation
$0\to \Ker\epsilon \To U \namedTo{\epsilon} M \to 0$ and a $R$-module
monomorphism $F\to G$ where $G$ is $R$-free, 
the sequence $\beta_1,\ldots, \beta_q \in U\backslash\Ker\epsilon$ 
together with $\varphi\in\Hom_R(U,R)$ 
satisfying the condition $(ii)$  in 
\keylemma 
is called a {\em \bseq} for the 
pair $(f:F\to G,M)$.
\end{Definition}

From \keylemma we immediagely 
have a characterization of long Bourbaki sequences.
\def\charBourbaki{Theorem~\ref{maincor}}
\begin{Theorem}
\label{maincor}
Let $r\in\ZZ$ be $r\geq 2$.
Let $R$ be a Noetherian  ring and $M$ be 
a finitely generated $R$-module.
Consider a $R$-module homomorphism
$f_1 : F_1\to M$ from a 
$R$-free module $F_1$.
Then, following are equivalent.
\begin{enumerate}
\item [$(i)$] We have a long Bourbaki sequence of length $r$
\[
 0 \To F_{r-1}\namedTo{f_{r-1}}\cdots\namedTo{f_2}F_1\namedTo{f_1} M
   \namedTo{\psi} I \To 0
\]
where  $I \subset R$ is an ideal and $F_i$ are $R$-free modules.
\item [$(ii)$]
There exists a \bseq $(\{\beta_i\}_i, \varphi)$ for 
$(\Ker f_1\hookrightarrow F_1, M)$ such that 
\begin{equation*}
 0 \To F_{r-1}\namedTo{f_{r-1}}\cdots\namedTo{f_3} F_2 \namedTo{f_2} \Ker f_1 \To 0	
\quad (exact)
\end{equation*}
\end{enumerate}
\end{Theorem}

\begin{Remark}
Notice that a b-sequence  in the case of length 2
Bourbaki sequence is not the same as a sequence of basic 
elements in the sense of \cite{Bour}. If we 
choose a suitable b-sequence $\{\beta_i\}$ 
under a suitabule condition, $\{\epsilon(\beta_i)\}$
can be a sequence of basic elements.
\end{Remark}

\subsection{Sygygies of Artinian Gorenstein rings}
A b-sequence has slightly more explicit description
for some class of modules 
over Gorenstein local rings.
Let $(R,\mm)$ be a Gorenstein local ring of dimension $n$.
We will denote the dual $\Hom_R(-, R)$ by $(-)^*$.  
Let $J_i\subset R\;(i=0,\ldots, d)$  $(d\leq n-r-1)$ be Gorenstein
ideals of grade $n$ and  set $M_i = R/J_i$.
Let $(F_\bullet^{(i)}, \partial_\bullet^{(i)})$ be 
a minimal $R$-free resolution of $M_i$.
By self-duality of the resolution we immediately have
\begin{Lemma}
\label{duality}
For all $i$, 
	$\Omega_i(M_i) \iso \Omega_{n-i+1}(M_i)^*
		       \iso \partial_{n-i+1}^*(\Omega_{n-i+1}(M)^*)$
\end{Lemma}
Now consider the module $M = \Dirsum_{i=0}^{d}\Omega_{i}(M_i)$.
By Lemma~\ref{duality} we have
\begin{Proposition}
\label{explicit}
Let $(\{\beta_i\}, \varphi)$ be a b-sequence for $M$. Then
$\varphi = \Dirsum_{i=0}^{d}a_i\circ \partial^{(i)}_i$
where $a_i \in \partial^*_{i}(\Omega_{i}(M_i)^*)$.
\end{Proposition}
This range of $a_i$ has more explicit description 
if $J_i =\mm =(x_1,\ldots,x_n)$ and $R = S = \poly{K}{x}{n}$.

For $I = \{i_1,\ldots, i_u\} \subset \{1,\ldots, n\}=[n]$,
we deonte by $e_I$ a base $e_{i_1}\wedge\cdots\wedge e_{i_u}$ of 
the Koszul complex $K_\bullet$ over $S$ of sequences $x_1,\ldots, x_n$.
A dual base to $e_I$ is denoted by $e_I^*$.
For $J, K \subset [n]$ with $J\cap K =\emptyset$ we
define $\sigma(J, K) = (-1)^i$ where $i = \sharp\{ (j,k)\in J\times K
\mid j > k\}$. Then we have $x_J\wedge x_K = \sigma(J,K)x_{J\cup K}$.
\begin{Corollary}
\label{moreexplicit}
Let $M_i = K (= R/\mm)$ for all $i$. Then $a_i$ in
Proposition~\ref{explicit} is an elemet from
\begin{equation*}
\left\langle 
\sum_{k=1}^{i}(-1)^{k+1} 
                   \sigma(J\backslash \{j_k\}, [n]-(J\backslash\{j_k\}))
                    x_{j_k} e^*_{[n]-(J\backslash\{j_k\})}
: J = \{j_1,\ldots, j_i\}\subset [n]
\right\rangle
\end{equation*}
\end{Corollary}

\section{Approximation of generalized Cohen-Macaulay ideals}

\def\unicond{homological condition}

\subsection{Approximation modules}
Let $(R,\mm)$ be a Gorenstein local ring and consider a generalized
Cohen-Macaulay ideal $I \subset R$ of codimension~$r$ ($r\geq 2$) such that
$\Coh{i}{R/I} = M_i$ for $i=0,\ldots, n-r-1$ where
$M_i$ are finite length $R$-modules.    
Then we have the following result, which is an immediate 
consequence from  Lemma~1.3 \cite{Ama1} and Theorem~1.1 \cite{HT2}.

\def\amahertak{Proposition~\ref{ama-her-tak}\ }
\begin{Proposition}
\label{ama-her-tak}
For a generalized Cohen-Macaulay ideal $I\subset R$ of codimension~$r
(\geq 2)$
there exists a maximal generalized Cohen-Macaulay module $M$ 
fitting into a Bourbaki sequence
\begin{equation*}
	0 \To F_{r-1}\To\cdots\To F_{1} \To M \To I \To I 
\end{equation*}
such that 
$\Coh{i}{M} \iso \Coh{i}{I}$ for $i\leq  n-r$, $\Coh{n-r+1}{M}=0$.
Moreover, 
if we assume the \unicond\ $(\ref{uniquenesscond})$,
then
$M$ is unique up to $R$-free direct summands.
\end{Proposition}

Notice that in \amahertak, the ideal $I$ is approximated 
by the module $M$ in a similar sense to Auslander-Buchweitz
(see \cite{HT2} for detail). We will call the maximal Cohen-Macaulay
module $M$ (or long Bourbaki sequence) an {\em approximation module}
(or {\em approximation sequence}).

More specific result can be obtained when we consider 
a special class of ideals.

\def\unique{Proposition~\ref{uniqueness}}
\begin{Proposition}
\label{uniqueness}
Let $(R,\mm)$ be a regular local ring and let 
$I\subset R$ be an ideal of codimension~$r$ $(\geq 2)$.
Assume that we have an approximation sequence
\begin{equation*}
0 \To F_{r-1}\To\cdots\To F_{1} \To \Omega_{t+1}(N)\dirsum H \To I \To 0
\end{equation*}
for some $R$-free modules $H$, $F_{1}, \ldots, F_{r-1}$ and 
a finite length $R$-module $N$. Then, we have 
\begin{equation}
\label{ss-cond}
	\Coh{i}{R/I} = 
	\left\{
		\begin{array}{ll}
		N  & i=t \\
	        0  & i< n-r, i\ne t\\
		\end{array}
	\right.
\end{equation}
Also the converse holds if
\begin{enumerate}
\item [(i)] $r =2$, or
\item [(ii)] $r\geq 3$ and we assume the \unicond\ $(\ref{uniquenesscond})$
for the approximation module $M$ of $I\subset R$.
\end{enumerate}
\end{Proposition}
\begin{proof}
The initial part of the proposition is clear. For the converse,
the case $r=2$ is Proposition~3.1 \cite{HT2}. 
Now assume that $I\subset R$ satisfies (\ref{ss-cond}) and 
the condition (ii). Then by \amahertak
we have an approximation sequence
\[
   0 \To F_{r-1}\To\cdots\To F_1\To X \To I \To 0
\]
with 
\[
   \Coh{i}{X}\iso \Coh{i}{I} \quad (i\le n-r),\qquad 
  \Coh{i}{X} = 0 \quad (n-r+1\le i \le n-1),
\]
and it remains to show that if $X$ is a $R$-module of maximal
dimension with the property that for $s := t+1$ with $0< s < n - r + 1$ one
has
\[
    \Coh{i}{X} \iso \left\{
			\begin{array}{ll}
0 & \mbox{for  $i<n$ and $i\ne s$} \\
N   & \mbox{for $i=s$}
			\end{array}
		\right.,
\]
then $X\iso \Omega_{s}(N)\dirsum H$ with some $R$-free module $H$.
But this fact is already proved in the proof of Proposition~3.1
\cite{HT2}.
\end{proof}

As proved in \amahertak and \unique, the \unicond\ (\ref{uniquenesscond})
assures the uniqueness of approximation
modules $M$. If we do not assume this condition, we have a large 
varieties of $M$ even in cohomologically very simple cases. For example,

\def\singlespot{Proposition~\ref{singlespot}}
\begin{Proposition}
\label{singlespot}
Let $r\geq 3$ and $0\leq t\leq n-r-1$ be integers.
Let $S = \poly{K}{x}{n}$ be a polynomial ring over a field 
$K$ and $\mm = (x_1,\ldots, x_n)$. 
Let $M$ be a maximal generalized CM module 
over $S$ with $\depth M = t+1$.
Consider a minimal $S$-free resolution of $M$:
\[
 F_\bullet : 
	0\To F_{n-t-1} \namedTo{\varphi_{n-t-1}}
	     F_{n-t-2} \namedTo{\varphi_{n-t-2}}
	         \cdots\namedTo{\varphi_2}
	     F_1       \namedTo{\varphi_1} 
             F_0       \namedTo{\varphi_0}
               M       \To 0.
\]
Also let $N$ be a non-zero finite length module over $S$.
Then the following are equivalent.
\begin{enumerate}
 \item [$(i)$]
  For any $l\in\ZZ$ such that $n-r+2 \leq l \leq n-1$,
  we have
\[ \Coh{i}{M} =
	\left\{
	\begin{array}{ll}
	   K   & \mbox{if } i=t+1 \\
	   N   & \mbox{if } i=l \\
           0   & \mbox{if } i<n, i\ne t+1, l\\
	\end{array}
      \right.
\]
\item [$(ii)$]
$\Omega_{n-l}(M) \iso F_{n-l}/E_{n+t+2-l}$, and
$N^\vee \iso \Omega_{n-l}(M)^*/\Im\varphi_{n-l}^*$.
\end{enumerate}
where we denote $\Omega_{s}(K)$ simply by $E_s$, and 
we define $(-)^* = \Hom_S(-,S(-n))$ and $(-)^{\vee} = \Hom_S(-,K)$.
\end{Proposition}
%
%
\begin{proof}
We first prove $(i)$ to $(ii)$.
By taking the dual of $F_\bullet$, we have 
\[
 0 \To F_0^*      \namedTo{\varphi_1^*} 
       F_1^*      \namedTo{\varphi_2^*}
       F_2^*      \namedTo{\varphi_3^*}
       \cdots     \namedTo{\varphi_{n-t-2}^*} 
       F_{n-t-2}^*\namedTo{\varphi_{n-t-1}^*}
       F_{n-t-1}^*\To 0.
\]
Then by local duality the  $j$th cohomology of this complex is
\[
   \Ext_S^j(M, S(-n)) \iso \Coh{n-j}{M}^\vee
   =   \left\{
	\begin{array}{ll}
		K      & \mbox{if $j=n-t-1$} \\
                N^\vee & \mbox{if $j=n-l$} \\
                0      & \mbox{if $j\ne n-t-1, n-l$}\\

	\end{array}
\right.
\] 
for $j\geq 1$. Thus 
\[
0\To \Im\varphi_{n+1-l}^* \To F_{n+1-l}^* \To\cdots\To F_{n-t-1}^*\To K \To 0
\]
is exact and $F_{n+1-l}^*$ to $F_{n-t-1}^*$ part is a begining of a minimal
free resolution of $K$, which is isomorphic to the 
corresponding  begining of the Koszul
complex $(K_\bullet,\partial_\bullet)$ of the sequence $x_1,\ldots,
x_n$. Namely,
\[
     F_{n-t-1}^* \iso K_0,\quad\ldots\quad,
     F_{n+1-l}^* \iso K_{l-t-2} \qquad\mbox{and}\quad 
     \Im\varphi_{n+1-l}^* \iso E_{l-t-1}.
\]
On the other hand, we have
$N^\vee \iso \Ker\varphi_{n+1-l}^*/\Im\varphi_{n-l}^*$ 
and 
$E_{l-t-1}\iso \Im\varphi_{n+1-l}^* \iso  F_{n-l}^*/\Ker\varphi_{n+1-l}^*$.
Now set $U := \Coker\varphi_{n-l}^* = F_{n-l}^*/\Im\varphi_{n-l}^*$. Then
\[
U/N^\vee \iso (F_{n-l}^*/\Im\varphi_{n-l}^*)/(\Ker\varphi_{n+1-l}^*/\Im\varphi_{n-l}^*)
 \iso  F_{n-l}^*/\Ker\varphi_{n+1-l}^* = E_{l-t-1}.
\]
Thus  we have
\[
      0 \To N^\vee \To U \To E_{l-t-1} \To 0
\]
Taking the dual, we have
\[
    0\To E_{l-t-1}^* \To U^* \To (N^\vee)^*.
\]
Since $N$ has finite length, $N^\vee$ has also finite length by Matlis
duality, so that $(N^\vee)^* =0$. 
Also $E_{l-t-1}^* \iso E_{n+t+2-l}$
by selfduality of Koszul complex. Thus we have $U^* \iso E_{n+t+2-l}$. 
Then by dualizing  the exact sequence
\begin{equation}
\label{prop:new1-seq1}
  F_{n-1-l}^* \namedTo{\varphi_{n-l}^*}F_{n-l}^* \To U \To 0
\end{equation}
we have
\[
	0\To E_{n+t+2-l} \To F_{n-l} \namedTo{\varphi_{n-l}} F_{n-1-l} 
         \To \Omega_{n-1-l}(M)\To 0.
\]
This proves the first condition of $(ii)$. 
Now from the short exact sequence
\[
    0 \To \Omega_{n-l}(M) \To F_{n-1-l} 
      \namedTo{\varphi_{n-1-l}} 
      \Omega_{n-1-l}(M) \To 0
\]
we have the long exact sequence 
\[
   0 \To \Omega_{n-1-l}(M)^* 
     \namedTo{\varphi^*_{n-1-l}}
         F_{n-1-l}^* \To 
         \Omega_{n-l}(M)^* \To
         N^\vee \To 0
\]
since we have $\Ext_S^{n-l}(M,S(-n))\iso \Coh{l}{M}^\vee = N^\vee$ by local
duality. Notice that we have $\Omega_{n-l}(M)^*\iso \Ker\varphi_{n+1-l}^*$
from the short exact sequence 
$F_{n+1-l}\namedTo{\varphi_{n+1-l}} F_{n-l}\To \Omega_{n-l}(M)\To 0$.
This proves the second condition in $(ii)$.

Next we prove $(ii)$ to $(i)$. By $(ii)(a)$ we have a $S$-free resolution
of $M$:
\[
   0\To K_n\namedTo{\partial_n}\cdots\namedTo{\partial_{n+t+3-l}}
          K_{n+t+2-l}\namedTo{\partial_{n+t+2-l}}
        F_{n-l}\namedTo{\varphi_{n-l}}\cdots\namedTo{\varphi_1}
        F_0 \namedTo{\varphi_0} M \To 0
\]
where $F_i$ are $S$-free modules. By taking the dual, we have
the complex
\[
    0\To 
      M^*    \namedTo{\varphi_0^*}
      F_0^* \namedTo{\varphi_1^*}\cdots \To
      F_{n-l}^* \namedTo{\partial_{n+t+2-l}^*}
	K_{n+t+2-l}^* \namedTo{\partial_{n+t+3-l}^*}\cdots\namedTo{\partial_n^*}
       K_n^* \To 0
\]
Then by local duality and selfduality of Koszul
complex we compute 
\[
    \Coh{i}{M} \iso \Ext_S^{n-i}(M,S(-n))^\vee
                =\left\{
			\begin{array}{ll}
			    K   & \mbox{if $i = t+1$} \\
                            0   & \mbox{if $i\leq l-1$, $i\ne t+1$}
			\end{array}
                \right.
\]
Now by dualizing the exact sequence
\[
    K_{n+t+2-l}\namedTo{\partial_{n+t+2-l}}
    F_{n-l} \namedTo{\varphi_{n-l}} \Ker\varphi_{n-1-l} \To 0
\]
we have 
\[
   0\To (\Ker\varphi_{n-1-l})^* \To F_{n-l}^* \namedTo{\partial_{n+t+2-l}^*} 
   K_{n+t+2-l}^*,
\]
so that we have $\Omega_{n-l}(M)^* =(\Ker\varphi_{n-1-l})^*\iso 
\Ker(\partial_{n+t+2-l}^*)$.
Then with the second condition in  $(ii)$ we compute
\begin{eqnarray*}
\Coh{l}{M}^\vee 
      &\iso &\Ext_S^{n-l}(M, S(-n)) 
        = \Ker \partial_{n+t+2-l}^* / \Im \varphi_{n-l}^*
\\
      &\iso & \Omega_{n-l}(M)^* / \Im \varphi_{n-l}^*
\\
      &  = & N^\vee
\end{eqnarray*}
as required.
\end{proof}

A typical class of the modules that do not satisfy the \unicond\ 
(\ref{uniquenesscond}) is 
$\Dirsum_{i=0}^{n-r-1}\Omega_{i+1}(M_i) \dirsum
\Dirsum_{i=n-r+2}^{n-1}\Omega_{i}(N_i)$ for finite 
length modules $M_i$ and $N_i$, which we will consider 
in the next subsection.

\subsection{Approximation sequences of non-trivial type}
In this subsection we assume $(R,\mm)$ to be regular local.
In the proof of \amahertak we construct approximation modules $M$
in a homological method, which always entails the \unicond\ 
(\ref{uniquenesscond}). See \cite{HT2} and \cite{Ama1}.
Now we are interested in the following problem: how can we construct
apprximation  sequences as in \amahertak that 
do not satisfy the \unicond\ (\ref{uniquenesscond})?
The simplest answer to this question is to make  the direct sum
of an approximation sequence as in \amahertak and 
the following exact sequences:
\begin{eqnarray*}
    0\To G^{(i)}_n \To \cdots \To G^{(i)}_{i} \To \Omega_{i}(N_i) \To 0	
\qquad (i=n-r+2,\ldots, n-1)
\end{eqnarray*}
where $N_i$ are any finite length $R$-modules and $G^{(i)}_\bullet$ are 
minimal $R$-free resolutions of $N_i$.
Then we have a Bourbaki sequence with the approximation module 
$M'= M\dirsum \Dirsum_{i=n-r+2}^{n-1}\Omega_{i}(N_i)$
and the map from $M'$ to the ideal $I$ is trivial on 
$\Dirsum_{i=n-r+2}^{n-1}\Omega_{i}(N_i)$ part. 
We will call this  an approximation 
sequence  of {\em trivial type}.

Now we will consider approximation sequences of non-trivial type.
Let $r\in\ZZ$ be $r\geq 2$ and $n\geq r+1$.
Consider a long Bourbaki sequence of length $r$
\begin{equation}
\label{bourbakiseq0}
    0 \To F_{r-1}\To\cdots\To F_2\To F_1\namedTo{g}
     M\dirsum N \namedTo{\phi} I \To 0
\end{equation}
where $I \subset R$ is a generalized Cohen-Macaulay ideal of codimension~$r$,
$F_i$  are $R$-free modules, and 
$M = \Dirsum_{i=0}^{n-r-1}\Omega_{i+1}(M_i)$
and $N = \Dirsum_{i= n-r+2}^{n-1}\Omega_{i}(N_i)$.
From this sequence, we construct the following
diagram, where $U_\bullet\namedTo{\epsilon} M$ and 
               $V_\bullet\namedTo{\eta} N$ are 
minimal free resolutions of $M$ and $N$,
and the third row is the mapping cone
$C(\alpha_\bullet)$ of a chain map $\alpha_\bullet$, which is a $R$-free
resolution of $I$. 
\begin{equation}
\label{mappingcone}
\begin{array}{cccccccl}
         &     &           &            &         &             &      0        &        \\
         &     &           &            &         &             &    \da        &        \\
  \cdots & \To &  F_2      &\To         & F_1     & \namedTo{g} & \Ker\phi      & \to 0  \\
         &     &\alpha_2\da&            &\alpha_1\da&           &  \da          & \      \\
  \cdots &\To  & U_{1}\dirsum V_1
                           & \To        & U_{0}\dirsum V_{0} 
                                                   & \namedTo{\epsilon\dirsum\eta} 
                                                                 &M\dirsum{N}   & \to 0\\
         &      &          &            &         &             &\hspace*{0.3cm}\da\phi  
                                                                                 &  \\
  \cdots &\To   & U_1\dirsum V_1\dirsum F_1 
                            &\To        & U_0\dirsum V_{0}
                                                   &  \To       &I              &\to 0 \\
          &     &           &           &         &             &  \da          &  \\
          &     &           &           &         &             &           0   &  \\
\end{array}
\end{equation}
Let $p_1 : U_0\dirsum V_{0} \to U_0$ and $p_2 : U_0\dirsum V_{0}\to V_{0}$ 
be the first and the second projections. 
From the diagram (\ref{mappingcone})  we know that 
    $\Ker\phi = \Im g = ((\epsilon\dirsum \eta)\circ \alpha_1)(F_1)$
and then by considering the ranks of the modules in the short exact sequence
\begin{equation}
\label{kerVarphi}
   0 \To \Ker\phi
     \To  M\dirsum {N}
     \namedTo{\phi} I
     \To 0
\end{equation}
we have
\begin{equation}
\label{rankcond1} 
    \rank(\Ker\varphi) = \rank(M) + \rank(N) - 1.
\end{equation}
On the other hand, we have 
\begin{equation}
\label{inclusion}
(\epsilon\circ p_1\circ\alpha_1)(F_1)
         \dirsum (\eta\circ p_2 \circ \alpha_1)(F_1)
  \supset
((\epsilon\dirsum\eta)\circ\alpha_1)(F_1)
 = \Ker\phi.
\end{equation}
Thus we have
\begin{eqnarray}
\label{rankcond2-1}
 \rank {\cal I}_{U_0} +\rank{\cal I}_{V_0} &\geq& \rank(M) + \rank(N) - 1 \\
\label{rankcond2-2}
 \rank(M) &\geq & \rank {\cal I}_{U_0}\\
\label{rankcond2-3}
 \rank(N) &\geq & \rank {\cal I}_{V_{0}} 
\end{eqnarray}
where  ${\cal I}_{U_0} := (\epsilon\circ p_1\circ \alpha_1)(F_1) (\subseteq M)$ and 
       ${\cal I}_{V_0} := (\eta    \circ p_2\circ \alpha_1)(F_1) (\subseteq N)$.
From this we know that 
$(\rank {\cal I}_{U_0}, \rank {\cal I}_{V_{0}})
= (\rank(M), \rank(N)), (\rank(M)-1, \rank(N))$, or 
$(\rank(M), \rank(N)-1)$.
Under this situation, we have
\begin{Lemma}
\label{chainmapBETA}
Following are equivalent.
\begin{enumerate}
\item [$(i)$]  the approximation sequence $(\ref{bourbakiseq0})$ is 
               of non-trivial type
\item [$(ii)$] For any free basis $\{m_i\}_i$ of $F_1$
there exists an index $i$ such that
$\alpha_1(m_i) \notin U_0$ and $\alpha_1(m_i)\notin V_{0}$.
\end{enumerate}
\end{Lemma}
\begin{proof}
We will prove $(i)$ to $(ii)$. 
We assume that for all $i$ we have either
$\alpha_1(m_i)\in U_{0}$ or $\alpha_1(m_i)\in V_{0}$, and 
will deduce  a contradiction.
First of all, we have equality in (\ref{inclusion}), and then
from (\ref{rankcond2-1}) we have 
\[
    \rank(M) + \rank(N) - 1 = \rank {\cal I}_{U_0} +  \rank {\cal I}_{V_{0}}.
\]
Thus, we have 
$(\rank {\cal I}_{U_{0}}, \rank {\cal I}_{V_{0}}) = 
(\rank(M)-1, \rank(N))$ or $(\rank(M), \rank(N) -1)$.
Also, since
$\Ker\phi = {\cal I}_{U_{0}}\dirsum {\cal I}_{V_{0}}$, we have by (\ref{kerVarphi})
\begin{equation}
\label{prop:I(c)}
I \iso (M/{\cal I}_{U_0})\dirsum (N/{\cal I}_{V_{0}})
\end{equation}
\begin{description}
\item [case $(\rank {\cal I}_{U_0}, \rank {\cal I}_{V_{0}}) = (\rank(M)-1, \rank(N)$]
Since we have
$\rank N/{\cal I}_{V_{0}} =\rank(N) - \rank {\cal I}_{V_{0}} = 0$,
$N/{\cal I}_{V_{0}}$ is 0 or a torsion-module. But since $I$ is torsion
free, we must have $N = {\cal I}_{V_{0}}$ by (\ref{prop:I(c)}). Thus 
$\Ker\phi = {\cal I}_{U_0}\dirsum N$ and then the Bourbaki sequence 
(\ref{bourbakiseq0}) must be of trivial-type
\[
     \cdots\To  F'_1 \dirsum V_{1} 
           \To  F'_0 \dirsum V_{0}
           \To  M \dirsum{N} 
            \namedTo{\phi} I\To 0
\]
where $F'_\bullet$ a $S$-free resolution 
of ${\cal I}_{U_0}$, a contradiction.

\item [case $(\rank {\cal I}_{U_0}, \rank {\cal I}_{V_{0}}) 
= (\rank(M), \rank(N) -1)$]
In this case we have $\rank M/{\cal I}_{U_0} =0$.
Since $M/{\cal I}_{U_0}\subset I$ by (\ref{prop:I(c)}) and $I$ is 
torsion-free, we must have $M/{\cal I}_{U_0} =0$. 
Thus $\Ker\phi = M\dirsum {\cal I}_{V_{0}}$ and 
the Bourbaki sequence (\ref{bourbakiseq0}) is obtained by
combining 
\begin{eqnarray*}
  0 \To M \dirsum {\cal I}_{V_{0}} \To M\dirsum N \namedTo{\phi} I \to 0
\end{eqnarray*}
with the minimal $R$-free resolution $U_\bullet$ of $M$
and a minimal $R$-free resolution of  ${\cal I}_{V_{0}}$.

Now since $M$ must satisfy $\Coh{i}{M} = \Coh{i}{I} \iso \Coh{i-1}{R/I}$ for 
$i = 1,\ldots, n-r$ and $\depth R/I \leq n-r-1$, we must have 
$\depth M \leq n-r$. Thus by Auslander-Buchsbaums formula, 
the length of $U_\bullet$ must be $\geq r$, which 
exceeds the length of our Bourbaki sequence, a contradiction.
\end{description}
Now we show $(ii)$ to $(i)$. 
Assume that  (\ref{bourbakiseq0}) is of trivial type. Then
we must have $\alpha_1(F_1) = p_1(\alpha_1(F_1))\dirsum p_2(\alpha_1(F_1))$. From
this we immediately obtain the required result.
\end{proof}

From Lemma~\ref{chainmapBETA}, we immediately have
\def\non-trivial-approx{Theorem~\ref{theorem:non-trivialcond}}
\begin{Theorem}
\label{theorem:non-trivialcond}
The appriximation sequence  (\ref{bourbakiseq0})
is of non-trivial type if and only if
\begin{enumerate}
\item [$(i)$] there exists a b-sequence 
	$\{\beta_i\}_i (\subset U_0\dirsum V_{0})$
              and
\item [$(ii)$] the submodule $N:= \langle \{\beta_i\}_i\rangle$
              of $U_0\dirsum V_{0}$ cannot be decomposed in the form
              of 
	      $N = A\dirsum B$ for some  $(0\ne) A\subset U_0$ and 
              $(0\ne )B \subset V_{0}$
\end{enumerate}
\end{Theorem}

\section{Some Applications  in Codimension~$3$}

\subsection{b-sequences for $E_{t+1}$ and $E_{t+1}\dirsum E_{n-1}(d)\;
 (d\in\ZZ)$}
As an application of our theory, we will consider a special case.
Let $S = \poly{K}{x}{n}$ and $\mm = (x_1,\ldots, x_n)$.
We consider the standard grading with $\deg(x_i)=1$ for all $i$.
Also, in the following, the dual $(-)^*$ always denotes 
$\Hom_S(-, S(-n))$.
We now consider the graded 
approximation module $M = E_{t+1}$ and $E_{t+1}\dirsum E_{n-1}(d)$,
for arbitrarily $d\in\ZZ$.

First of all, by \keylemma and Corollary~\ref{moreexplicit} we 
have the following.
\def\explicit1{Corollary~\ref{explicit1}\ }
\begin{Corollary}
\label{explicit1}
Following are equivalent.
\begin{enumerate}
\item [$(i)$] We have a length $3$ Bourbaki sequence
\[
 0 \To F \namedTo{f} G \namedTo{g}  E_{t+1}\dirsum E_{n-1}(d)
   \To I(c) \To 0 \quad (exact)
\]
where  $I \subset S$ is a graded ideal and $F$ and $G$ are 
finitely generated $S$-free modules.
\item [$(ii)$]
$\rank F = \rank G - n + 2 - \binom{n-1}{t}$ and 
we have a b-sequence 
$(\{\beta_i\}_i, \varphi)$ 
for $(f, E_{t+1}\dirsum E_{n-1}(d))$ where 
$\beta_i\in K_{t+1}\dirsum K_{n-1}(d)\backslash E_{t+2}\dirsum E_n(d)$
and $\varphi = (a,b) \in {\cal A}\times{\cal B}$, with 
\begin{eqnarray*}
{\cal A}& =& \left\langle
		\sum_{j=1}^{n-t}
			(-1)^{j+1}\sigma(L\backslash\{i_j\}, 
([n]\backslash L)\cup\{i_j\})
				x_{i_j} e^*_{([n]\backslash L)\cup\{i_j\}}
		\mid L = \{i_1,\ldots,i_{n-t}\}\subset [n]
            \right\rangle \\
{\cal B} &= &\langle
		(-1)^ix_j e^*_{[n]\backslash\{i\}} -
		(-1)^jx_i e^*_{[n]\backslash\{j\}}
		\mid 1\leq i < j \leq n
            \rangle,
\end{eqnarray*}
and thus $\varphi :  K_{t+1}\dirsum K_{n-1}(d)\to S(-n)$ is 
      a degree '$n+c$' homomorphism.
\end{enumerate}
In this case, we have $I = \varphi(K_{t+1}\dirsum K_{n-1}(d))(-c)$
\end{Corollary}

We also consider the case of $M = E_{t+1}$.
\begin{Corollary}
\label{explicit2}
Following are equivalent.
\begin{enumerate}
\item [$(i)$] We have a length $r$ $(\geq 3)$
Bourbaki sequence
\[
 0 \To F_{r-1}\To\cdots\To F_2\To F_1
   \namedTo{f_1}  E_{t+1}
   \To I(c) \To 0 \quad (exact)
\]
where  $I\subset S$ is a graded ideal, and 
$F_i$ are finitely generated $S$-free modules.
\item [$(ii)$]
We have 
$\rank \Ker f_1 = \rank F_1 + 1 - \binom{n-1}{t}$ and 
a b-sequence  $(\{\beta_i\}_i \subset 
K_{t+1}\backslash E_{t+2}, \varphi\in {\cal A})$
for $(\Ker f_1\hookrightarrow F_1, E_{t+1})$
where 
\[
{\cal A} = \left\langle
		\sum_{j=1}^{n-t}
			(-1)^{j+1}\sigma(L\backslash\{i_j\},
 ([n]\backslash L)\cup\{i_j\})
				x_{i_j} e^*_{([n]\backslash L)\cup\{i_j\}}
		\mid L=\{i_1,\ldots,i_{n-t}\} \subset [n] \right\rangle
\]
and thus $\varphi : K_{t+1} \to S(-n)$ defines a degree '$n+c$' homomorphism.
In this case, we have $I = \varphi(K_{t+1})(n-c)$.
\end{enumerate}
\end{Corollary}

A small application of this explicit formula is 
\begin{Corollary}
\label{cor:depthzero}
There is no graded ideal $I\subset S$ of codimension $r$ $(\geq 2)$
of $\depth(S/I) = 0$ such that local cohomology is 
trivial except  $\Coh{0}{S/I} = K(c)$ (for some $c\in\ZZ$)
and having a length $r$  approximation sequence with 
approximation module $E_1$.
\end{Corollary}
\begin{proof}
Assume that there exists an ideal $I \subset S$ 
such that $\Coh{0}{S/I} = K(c)$ and $\Coh{i}{S/I}=0$ 
$(0<i<n-r)$ having the following 
approximation seqence:
\[
   0 \To F_{r-1} \To\cdots\To F_1 \To E_1 \namedTo{\varphi} I(c) \To 0.
\]
Then by Corollary~\ref{explicit2}, there exists a b-sequence
$(\{\beta_i\}_i, \varphi)$ with 
$\beta_i \in K_1\backslash E_2$ and 
$(0\ne) \varphi \in {\cal A}$.
Then $\langle\{\beta_i\}_i\rangle + E_2 
= \Ker(\varphi : K_1\to S(-n))$. 
Since $\varphi$ is a non-zero element from
${\cal A} =  \langle x_1e^*_1 +\cdots +  x_ne^*_n\rangle$
we must have $\{\beta_i\}_i \subset  E_2(\ne E_1)$, 
a contradiction.
\end{proof}

\begin{Remark}
By \unique\ it is assured that an ideal as in 
Corollary~\ref{cor:depthzero} has a length $r$
approximation sequence with approximation module
$E_1\dirsum H$, with non-trivual $S$-free module $H$. 
We will show later that 
there exists an approximation sequence with 
approximation module $E_1\dirsum E_{n-1}$, due to 
Corollary~\ref{explicit1}. See Example~\ref{example3}.
\end{Remark}

\subsection{Numerical condition for codimension~3}
Now we consider in par ticular the case of \explicit1. 
The existence of approximation sequences as in \explicit1
only implies that $\codim I \leq 3$.
To assure that the codimension is exactly $3$, we need
additional condition. We have 
\def\numerical{Proposition~\ref{theorem:numerical}}
\begin{Proposition}
\label{theorem:numerical}
Let $n\geq 4$ and $t\leq n-4$. Assume that we have the following long Bourbaki sequence
\begin{equation}
\label{codim3bseq}
0\To \Dirsum_{i=1}^p S(-a_i)\To \Dirsum_{i=1}^q S(-b_i)\To E_{t+1}\dirsum E_{n-1}(d)\To I(c)\To 0
\end{equation}
with $I \subset S$ a graded ideal and $c\in\ZZ$. Then we have $\codim
I \leq 3$ and the equality holds if and only if
\begin{enumerate}
\item $q = p + \binom{n-1}{t} + n - 2$;
\item $\displaystyle{\sum_{i=1}^q }b_i - \displaystyle{\sum_{i=1}^p} a_i
 = n^2 - (2+d)n + c + d + \binom{n-2}{t-1} + \binom{n-1}{t}t$;  
\item and \begin{eqnarray*}
\displaystyle{\sum_{i=1}^q} b^2_i - \displaystyle{\sum_{i=1}^p} a^2_i
               &=& n^3 - (3+2d)n^2 + (d^2+4d+1)n - c^2 - d^2 \\
               & & + \binom{n-1}{t}(t+1)^2 -\binom{n-2}{t}(2t+1) - 2\binom{n-3}{t-1}
      \end{eqnarray*}
\end{enumerate}
\end{Proposition}
\begin{proof}
Now from the sequence (\ref{codim3bseq}),  we construct 
the mapping cone $C(\alpha_\bullet)$ in a similar way to (\ref{mappingcone}).
The cone gives a $S$-free resolution $F_\bullet$ of the residue ring $S/I$.
\[
F_\bullet: 0\to F_{n-t}\to\cdots\to F_1 \to F_0 \to S/I \to 0
\]
where 
\begin{eqnarray*}
F_0 & = & S \\
F_1 & = & K_{t+1}(-c)\dirsum K_{n-1}(d-c) = S(-t-1-c)^{\beta_1}\dirsum S(-n+1+d-c)^n  \\
F_2 & = & K_{t+2}(-c)\dirsum K_n(d-c)\dirsum G(-c) \\
    & = & S(-t-2-c)^{\beta_2}\dirsum S(-n+d-c)\dirsum \Dirsum_{i=1}^{q}S(-b_i-c)\\
F_3 & = & K_{t+3}(-c)\dirsum F(-c) = S(-t-3-c)^{\beta_3}\dirsum \Dirsum_{i=1}^{p}S(-a_i-c)\\
F_i & = & K_{t+i}(-c) =S(-t- i -c)^{\beta_{i}} 
\qquad  (4\leq i\leq n-t)\\
\mbox{with} && \beta_i = \binom{n}{t+i}\quad i=1,\ldots, n-t.
\end{eqnarray*}
Now we compute the Hilbert series $Hilb(S/I, \lambda)$ of $S/I$. We have 
\begin{equation}
\label{hilbertSeries}
    Hilb(S/I, \lambda) = \frac{ Q(\lambda)}{(1-\lambda)^n}
\end{equation}
with
\begin{eqnarray*}
\lefteqn{Q(\lambda) =  \sum_{i,j}(-1)^i\beta_{i,j}\lambda^j} \\
           &=& 1 - n \lambda^{n-1+c-d} + \lambda^{n+c-d} 
                + \sum_{i=1}^q\lambda^{b_i+c}
- \sum_{i=1}^p \lambda^{a_i+c}
+ (-1)^t\lambda^c\sum_{i=t+1}^n(-1)^i\binom{n}{i}\lambda^i
\end{eqnarray*}
where $\beta_{i,j}$ are as in $F_i = \Dirsum_{j}S(-j)^{\beta_{ij}}$, $(i=0,\ldots, n-t)$.
(see Lemma~4.1.13 \cite{BH}).
Since we have $\Coh{i}{S/I} = \Coh{i+1}{M}$ for $0\leq i \leq n-4$, 
we know that $\dim S/I \geq n-3$, i.e.,
$\codim I \leq 3$. To assure that $\codim I \geq 3$ we must have $Q(1)
= Q'(1) = Q''(1) = 0$ (see Corollary~4.1.14(a) \cite{BH}).

Now by straightforward computations, we have 
\begin{enumerate}
\item $Q(1) = 0$ holds for all $n,t,c$ and $p$
\item $Q'(1) = 0$ holds if and only if 
\[
   \sum_{i=1}^q b_i  - \sum_{i=1}^p a_i 
= n^2 -(2+d)n + c + d + \binom{n-2}{t-1}
                         + \binom{n-1}{t}t.
\]
\item $Q''(1)=0$ holds if and only if 
\begin{eqnarray*}
\sum_{i=1}^q b_i^2  - \sum_{i=1}^p a_i^2
&=& n^3 - (3+2d)n^2  + (d^2+4d+1)n - c^2 - d^2  \\
& &
 + \binom{n-1}{t}(t+1)^2 
    - \binom{n-2}{t}(2t+1) - 2\binom{n-3}{t-1}
\end{eqnarray*}
\end{enumerate}
as required.
\end{proof}

\subsection{Examples}
Now we give a few concrete examples in codimension~3.

\begin{Example}[approximation module $E_{t+1}$ with $t=\depth S/I = 1$]
\label{example1}
We first give an application of Corollary~\ref{explicit2}. Namely,
a codimension~3 ideal $I$ with approximation module $E_{t+1}$. 
Let $t = 1$ and $n =6$. Then ${\cal A} = \{A_1,A_2,A_3,A_4,A_5,A_6\}
 = \partial_2^*(E_2^*) \subset K_2^*$ where 
\begin{eqnarray*}
A_1 & = &  x_1 e^*_{16} + x_2 e^*_{26} + x_3 e^*_{36} + x_4 e^*_{46} + x_5 e^*_{56}
\\
A_2 & = & -x_1 e^*_{15} - x_2 e^*_{25} - x_3 e^*_{35} - x_4 e^*_{45} + x_6 e^*_{56}
\\
A_3 & = &  x_1 e^*_{14} + x_2 e^*_{24} + x_3 e^*_{34} - x_5 e^*_{45} - x_6 e^*_{46}  
\\
A_4 & = & -x_1 e^*_{13} - x_2 e^*_{23} + x_4 e^*_{34} + x_5 e^*_{35} + x_6 e^*_{36}  
\\
A_5 & = &  x_1 e^*_{12} - x_3 e^*_{23} - x_4 e^*_{24} - x_5 e^*_{25} - x_6 e^*_{26}  
\\
A_6 & = &  x_2 e^*_{12} + x_3 e^*_{13} + x_4 e^*_{14} + x_5 e^*_{15} + x_6 e^*_{16}.
\end{eqnarray*}
We choose a b-sequence $(\{\beta_i\}_i, a)$ with 
$a \in {\cal A}$ and $\beta_i\in K_2\backslash E_3$ 
as follows:
\begin{eqnarray*}
a & = & x_6 A_1 -x_5 A_2 + x_4 A_3 \\
  & = & x_1x_4 e^*_{14} + x_1x_5 e^*_{15} +x_1x_6 e^*_{16} + x_2x_4 e^*_{24}
       + x_2x_5 e^*_{25}\\
  &   & + x_2x_6 e^*_{26} +x_3x_4 e^*_{34} + x_3x_5 e^*_{35}
       sy + x_3x_6 e^*_{36} \\
\beta_1 & = & e_{12},\quad \beta_2 = e_{13},\quad \beta_3 = e_{23},
                     \quad \beta_4 = e_{45},\quad \beta_5 = e_{46},
                     \quad \beta_6 = e_{56}
\end{eqnarray*}
Then we obtain the long Bourbaki sequence
\[
   0 \To S^2(-3) \namedTo{f} S^6(-2) \namedTo{g} E_2 \namedTo{\varphi} I \To 0.
\]
We can check that all the conditions in \numerical\ are satisfied 
so that we must have $\codim I = 3$.
If $\{n_i\}$ and $\{m_j\}$ are free bases of $S^2(-3)$ and 
$S^6(-2)$, we have 
\[
\begin{array}{cccc}
f : & S^2(-3) = S n_1 \dirsum S n_2 & \To &  S^6(-2) = S m_1 \dirsum\cdots\dirsum S m_6 \\
    &           n_1                 &\longmapsto  &  x_3 m_1 - x_2 m_2 + x_1 m_3 \\
    &           n_2                 &\longmapsto  &  x_6 m_4 - x_5 m_5 + x_4 m_6 \\
    &                               &             &                               \\
g : &  S^6(-2) = S m_1 \dirsum\cdots\dirsum S m_6 & \To  & E_2 \\
    &           m_i                 &\longmapsto  &  \partial_2(\beta_i) \qquad (i=1,\ldots,6) \\
\mbox{and} &                        &             &                               \\
\varphi &   E_2               & \To         & I  \\
        & \partial_2(e_{ij})  & \longmapsto & x_ix_j \\
        &                     &\mbox{for }& (i,j) = (1,2),(1,3),(2,3),(4,5),(4,6),(5,6)\\
        & \partial_2(e_{ij})  & \longmapsto & 0 \\
        &                     &\mbox{for }& (i,j) \ne (1,2),(1,3),(2,3),(4,5),(4,6),(5,6)\\
\end{array}
\]
and we obtain 
$I = (x_1x_4,x_1x_5, x_1x_6,x_2x_4,x_2x_5,x_2x_6,x_3x_4, x_3x_5, x_3x_6)
       = (x_1,x_2,x_3)(x_4,x_5,x_6)$.
\end{Example}

\begin{Example}[approximation module $E_{t+1}\dirsum E_{n-1}$ with $t= \depth S/I = 1$]
\label{example2}
We continue to consider the situation in Example~\ref{example1}. 
As an application of Corollary~\ref{explicit1}, we can see
that the same ideal fits into a long Bourbaki sequence with approximation
module $E_{t+1}\dirsum E_{n-1} = E_1 \dirsum E_5$.
In this case, we must also  consider 
${\cal B} = \{B_{ij} \mid 1\leq i <j \leq 6 \}
 = \partial_5^*(E_5^*) \subset K_5^*$ 
where $B_{ij} = (-1)^{i}x_je_{[6]\backslash i}^*
              - (-1)^{j}x_ie_{[6]\backslash j}^*$.
Then we set $a \in {\cal A}$ as in Example~\ref{example1} and 
\[
b = -x_1^2x_2x_4 B_{14} =x_1^2x_2x_4^2 e^*_{23456} + x_1^3x_2x_4 e^*_{12356}
\in {\cal B}.
\]
Also we set $\beta_1,\ldots, \beta_6$ to be the same as 
those in Example~\ref{example1} and 
$\beta_7 =x_1x_2x_4 e_{14} - e_{23456}$,
$\beta_8 =x_1^2x_2 e_{14} - e_{12356}$,
$\beta_9    = e_{13456}$,
$\beta_{10} = e_{12456}$,
$\beta_{11} = e_{12346}$,
$\beta_{12} = e_{12345}$.
We can check that $(\{\beta_i\}_i, \varphi = (a,b))$ is 
a b-sequence for $(N\hookrightarrow K, E_1\dirsum E_5)$
where $L$ is a free module of rank $12$ and 
$N$ is 
\begin{equation}
N = 
\left\langle
          \begin{array}{l}
            x_3m_1-x_2m_2+x_1m_3,\;
            x_6m_4-x_5m_5+x_4m_6,\\
          x_1m_7-x_4m_8+x_2m_9-x_3m_{10}-x_5m_{11}+x_6m_{12}
          \end{array}
         \right\rangle
\end{equation}
where $\{m_i\}_i$ is a free basis of $L$.
Notice that $\{\beta_i\}_{i=1}^{12}$ satisfies the condition of
Proposition~\ref{theorem:non-trivialcond}.

Then we have an approximation sequence of non-trivial type
\[
   0 \To S^2(-3)\dirsum S(-6) 
    \namedTo{f} 
         S^6(-2)\dirsum S^6(-5)
    \namedTo{g} 
         E_2\dirsum E_5
    \namedTo{\varphi} I \To 0
\]
where 
\[
\begin{array}{cccc}
f : & S^2(-3)\dirsum S(-6) = \Dirsum_{i=1}^3S n_i
                                    & \To         &  S^6(-2)\dirsum S^6(-5) = \Dirsum_{i=1}^{12}S m_i \\
    &           n_1                 &\longmapsto  &  x_3 m_1 - x_2 m_2 + x_1 m_3 \\
    &           n_2                 &\longmapsto  &  x_6 m_4 - x_5 m_5 + x_4 m_6 \\
    &           n_3                 &\longmapsto  &  x_1 m_7 - x_4 m_8 + x_2 m_9 \\
    &                               &             & - x_3 m_{10} -x_5 m_{11} + x_6 m_{12} \\
    &                               &             &                               \\
g : &   S^6(-2)\dirsum S^6(-5) = \Dirsum_{i=1}^{12}S m_i
                                    & \To  & E_2\dirsum E_5 \\
    &           m_i                 &\longmapsto  &  \bar\partial(\beta_i) \qquad (i=1,\ldots,12) \\
    &                               &             & \mbox{where }\bar\partial = \partial_2\dirsum\partial_5\\
\mbox{and} &                        &             &                               \\
\varphi &   E_2               & \To         & I  \\
        & \partial_2(e_{ij})  & \longmapsto & x_ix_j \\
        &                     &\mbox{for }& (i,j) = (1,2),(1,3),(2,3),(4,5),(4,6),(5,6)\\
        & \partial_2(e_{ij})  & \longmapsto & 0 \\
        &                     &\mbox{for }& (i,j) \ne (1,2),(1,3),(2,3),(4,5),(4,6),(5,6)\\
        & \partial_5(e_{23456}) & \longmapsto & x_1^2x_2x_4^2 \\
        & \partial_5(e_{12356}) & \longmapsto & x_1^3x_2x_4 \\
        & \partial_5(e_{ijklm}) & \longmapsto & 0\quad \mbox{otherwise} \\
\end{array}
\]
and the ideal $I$ is the same as that in Example~\ref{example1}.
We can also check that this sequence satisfies the numerical condition in Theorem~\ref{theorem:numerical}
\end{Example}

\begin{Example}[approximation module $E_{t+1}\dirsum E_{n-1}(d)$ with 
$t=\depth S/I = 0$]
\label{example3}
By Corollary~\ref{cor:depthzero}, we do not have a long Bourbaki
sequence with an approximation module $E_1$ and a codimension~$3$
generalized CM ideal $I$. However, there are long Bourbaki sequences
with approximation modules $E_1\dirsum E_5(d)$ for $d\in\ZZ$, which 
is an application of Corollary~\ref{explicit1}.
Let $k=1$ and $n=6$. Then,
we choose a b-sequence $(\{\beta_i\}_i,\varphi)$ as follows:
We set $\beta_i\in K_1\dirsum K_5(1)$ to be 
\begin{equation*}
\begin{array}{cclccl}
\beta_1 & = & - x_6e_{12345}  + x_5 e_{12346}, &   
\beta_2 & = &  x_6^5 e_3 - x_1^2 e_{13456},    \\
\beta_3 & = &  x_6^5 e_2 - x_1^2 e_{12456},    &
\beta_4 & = &  x_2^4x_5 e_2 - x_1 e_{12345},   \\
\beta_5 & = &  x_2^4x_6 e_2 - x_1^2 e_{12346}, &
\beta_6 & = & -x_6^4 e_{12346} + x_2^4 e_{12456}, \\
\beta_7 & = & e_{23456}, &
\beta_8 & = & e_{12356}. \\
\end{array}
\end{equation*}
Also let $\varphi = (a,b) \in {\cal A}\times{\cal B}$ be
\begin{eqnarray*}
a & = & x_1^3e_1^* + x_1^2x_2 e_2^* + x_1^2x_3 e_3^* + 
        x_1^2x_4 e_4^* + x_1^2x_5 e_5^* + x_1^2x_6 e_6^* \\
b & = & -x_2^5 B_{56} + x_6^5 B_{23}\\
  & = & x_2^5x_6e^*_{12346} + x_2^5x_5 e^*_{12345}
         +x_3x_6^5 e^*_{13456} + x_2x_6^5 e^*_{12456} \\
\end{eqnarray*}
where 
\begin{eqnarray*}
{\cal A} & =&  \langle
                   x_1e_1^* + \cdots + x_6 e_6^* 
               \rangle\\
{\cal B} &= & \langle B_{ij} = (-1)^ix_je_{[6]\backslash i}^*
- (-1)^jx_ie_{[6]\backslash j}^*  :
		  1\leq i < j \leq 6 
\rangle.
\end{eqnarray*}
Then we can check 
that $(\{\beta_i\}_i,\varphi)$ is a b-sequence for
$(F'\hookrightarrow G, E_1\dirsum E_5(1))$
where 
\begin{eqnarray*}
G & = & \langle m_1,\ldots, m_8\rangle
=S(-5)\dirsum S^4(-6)\dirsum S(-8)\dirsum S^2(-4) \\
F' & = & \left\langle 
\begin{array}{l}
-x_1^2 m_1 + x_6 m_4 - x_5 m_5,
                 x_2^4 m_3 - x_6^4 m_5 + x_1^2 m_6,\\
                -x_1^2m_1 -x_2 m_2 + x_3 m_3 - x_1^3 m_7 + x_1^2x_4 m_8
             \end{array}
         \right\rangle \\
\end{eqnarray*}
Also we know that the condition of Theorem~\ref{theorem:non-trivialcond} 
is satisfied.
Thus we have a non-trivial approximation sequence
\[
  0 \to F \namedTo{f} G \namedTo{g} E_1\dirsum E_5(1) \namedTo{\phi} I(2) \to 0
\]
where $g(m_i) = \beta_i$, $i=1,\ldots, 8$, and $F = S(-10)\dirsum S^2(-7)
= \langle u, v, w\rangle$ with 
$f(u) = x_2^4 m_3 - x_6^4 m_4 + x_1^2 m_6$,
$f(v) = -x_1^2 m_1 + x_6 m_4 - x_5m_5$
and 
$f(w) = -x_1^2 m_1 - x_2 m_2
+ x_3 m_3 - x_1^3 m_7 + x_1^2x_4 m_8$. The map $\phi$ is as follows:
$\phi(x_i) = x_ix_1^2$ $(i=1,\ldots, 6)$, $\phi(\partial_5(e_{12345})) = x_2^5x_5$,
$\phi(\partial_5(e_{12346})) = x_2^5x_6$, $\phi(\partial_5(e_{12356})) = 0$,
$\phi(\partial_5(e_{12456})) = x_2x_6^5$, $\phi(\partial_5(e_{13456})) 
   = x_3x_6^5$, and $\phi(\partial_5(e_{23456})) = 0$. 
The ideal is $I = \Im\varphi = x_1^2\mm + (x_2^5x_6, x_2^5x_5, x_3x_6^5, x_2x_6^5)$.
Finally we can check that this approximation sequence satisfies the numerical
condition of Theorem~\ref{theorem:numerical}, so that $\codim I = 3$.
\end{Example}


\end{document}